\documentclass[a4paper]{amsart}
\usepackage{color}
\usepackage[latin1]{inputenc}
\usepackage[T1]{fontenc}
\usepackage{amsfonts}
\usepackage{amssymb}
\usepackage{amsmath}
\usepackage{amsthm}
\usepackage{stmaryrd}
\usepackage{enumerate}
\usepackage{multirow}
\usepackage{graphicx}
\usepackage{setspace}
\usepackage{graphicx,type1cm,eso-pic,color}
\usepackage{pstricks}
\usepackage{float}
\newtheorem{theorem}{Theorem}[section]
\newtheorem{lemma}[theorem]{Lemma}
\newtheorem{proposition}[theorem]{Proposition}

\newtheorem{assumption}[theorem]{Assumption}
\newtheorem{remark}[theorem]{Remark}

\begin{document}
\setlength\arraycolsep{2pt}
\title[On a Universal Strictly Decreasing Nonparametric Estimator]{On a Universal Strictly Decreasing Nonparametric Estimator Applied to the Drift Function of a Recurrent Diffusion Process Estimation}
\author{Nicolas MARIE$^{\dag}$}
\address{$^{\dag}$Modal'X, Universit\'e Paris Nanterre, Nanterre, France}
\email{nmarie@parisnanterre.fr}
\date{}
\maketitle
%


%
\begin{abstract}
This paper deals with a copies-based continuously differentiable and strictly decreasing estimator of the drift function for stochastic differential equations defining recurrent diffusion processes. The first part of our paper deals with non-asymptotic $\mathbb L^1$-risk bounds and a bandwidths selection procedure for a universal monotone estimator. These results are tailor-made to our framework, and then applied to the estimation of the drift function of recurrent diffusion processes in the second part of the paper.
\end{abstract}
\tableofcontents
%


%
\section{Introduction}\label{section_introduction}
Consider the stochastic differential equation
\begin{equation}\label{main_equation}
X_t = x_0 +
\int_{0}^{t}a(X_s)ds +
\int_{0}^{t}\sigma(X_s)dW_s,
\quad t\in\mathbb R_+,
\end{equation}
where $x_0\in\mathbb R$, $W$ is a Brownian motion, $a,\sigma\in C^1(\mathbb R)$ and $a'$, $\sigma$ and $\sigma'$ are bounded. Under these conditions, Equation (\ref{main_equation}) has a unique (strong) solution $X$. Moreover, under the following Assumption \ref{assumption_drift_volatility} on $(a,\sigma)$, the scale density
\begin{displaymath}
\texttt s(\cdot) =\exp\left(-2\int_{\alpha}^{\cdot}\frac{a(y)}{\sigma(y)^2}dy\right)
\quad (\alpha = a^{-1}(0))
\end{displaymath}
satisfies
\begin{displaymath}
\int_{-\infty}^{\alpha}\texttt s(x)dx =
\int_{\alpha}^{\infty}\texttt s(x)dx =\infty,
\end{displaymath}
and then $X$ is a recurrent Markov process (see Gihman and Skorokhod \cite{GS72}, Theorem 1 p. 162).
%


%
\begin{assumption}\label{assumption_drift_volatility}
There exist $\mathfrak m_a,\mathfrak m_{\sigma} > 0$ such that, for every $x\in\mathbb R$, $a'(x)\leqslant -\mathfrak m_a$ and $|\sigma(x)|\geqslant\mathfrak m_{\sigma}$.
\end{assumption}
\noindent
Since, up to our knowledge, there is no reference about statistical inference for stochastic differential equations involving a (strictly) monotone drift function; our purpose is to investigate a transformation, taking into account both the regularity and the monotonicity conditions on $a$, of the copies-based Nadaraya-Watson estimator provided in Marie and Rosier \cite{MR23}.
\\
\\
{\bf Observation scheme:} Consider $N$ independent copies $X^1,\dots,X^N$ of $X_{|[0,T]}$ such that
\begin{equation}\label{observation_scheme}
X^i =\mathcal I(x_0,W^i),
\quad\forall i\in\{1,\dots,N\},
\end{equation}
where $N\in\mathbb N^*$, $T > 0$, $\mathcal I(\cdot)$ is the It\^o map for Equation (\ref{main_equation}), and $W^1,\dots,W^N$ are independent copies of $W_{|[0,T]}$. For instance, $X^1,\dots,X^N$ are appropriate to model the elimination processes of a drug by $N$ patients involved in a clinical trial (see Donnet et Samson \cite{DS13}). Moreover, under Assumption \ref{assumption_drift_volatility}, $X^1,\dots,X^N$ can be constructed from one observation of $X$ on $\mathbb R_+$. To that purpose, consider the stopping times $\tau_1,\dots,\tau_N$, recursively defined by $\tau_1 = 0$ and
\begin{displaymath}
\tau_i =\inf\{t >\tau_{i - 1} + T : X_t = x_0\},
\quad i = 2,\dots,N.
\end{displaymath}
For any $i\in\{1,\dots,N\}$, $\mathbb P(\tau_i <\infty) = 1$ because $X$ is a recurrent Markov process under Assumption \ref{assumption_drift_volatility}, and then
\begin{displaymath}
W^i := (W_{\tau_i + t} - W_{\tau_i})_{t\in [0,T]}
\quad {\rm and}\quad
X^i := (X_{\tau_i + t})_{t\in [0,T]}
\quad\textrm{are well-defined}.
\end{displaymath}
So, $W^1,\dots,W^N$ are independent Brownian motions by the strong Markov property, and $X^1,\dots,X^N$ are independent copies of $X_{|[0,T]}$ satisfying (\ref{observation_scheme}).
\\
\\
For stochastic differential equations with no monotonicity condition on the drift function, several nonparametric copies-based estimators of $a$ have been investigated in the literature. For instance, Comte and Genon-Catalot \cite{CGC20} (resp. Comte and Marie \cite{CM23}) deals with a projection least squares estimator of the drift function computed from independent (resp. correlated) continuous-time observations of a diffusion process, Denis et al. \cite{DDM21} with such an estimator computed from independent discrete-time observations, and Marie and Rosier \cite{MR23} with continuous-time and discrete-time versions of a copies-based Nadaraya-Watson estimator of $a$. On such estimators for interacting particle systems, the reader can refer to Della Maestra and Hoffmann \cite{DMH22}, while Amorino et al. \cite{ABPPZ25} investigates the estimator of the interaction function for a class of McKean-Vlasov stochastic differential equations when both $T$ and $N$ goes to infinity. Finally, even if our paper focuses on the copies-based observation scheme, let us mention some references dealing with nonparametric estimators of the drift function computed from one long-time observation of the ergodic stationary solution of Equation (\ref{main_equation}): Comte, Genon-Catalot and Rozenholc \cite{CGCR07} (resp. Aeckerle-Willems and Strauch \cite{AWS22}) on a discrete-time (resp. continuous-time) projection least squares (resp. Nadaraya-Watson) estimator of $a$.
\\
\\
{\bf Estimators:} Consider the Nadaraya-Watson estimator
\begin{displaymath}
\widehat a_{\eta}(\cdot) =
\frac{\widehat{af}_{\eta}(\cdot)}{\widehat f_{\eta}(\cdot)}
\mathbf 1_{\widehat f_{\eta}(\cdot) >\frac{\mathfrak m}{2}}
\end{displaymath}
with $\mathfrak m > 0$ defined later,
\begin{displaymath}
\widehat{af}_{\eta}(\cdot) =
\frac{1}{NT_0}\sum_{i = 1}^{N}\int_{t_0}^{T}K_{\eta}(X_{s}^{i} -\cdot)dX_{s}^{i}
\quad {\rm and}\quad
\widehat f_{\eta}(\cdot) =
\frac{1}{NT_0}\sum_{i = 1}^{N}\int_{t_0}^{T}K_{\eta}(X_{s}^{i} -\cdot)ds,
\end{displaymath}
where $T_0 = T - t_0$, $t_0\in (0,T)$,
\begin{displaymath}
K_{\eta}(\cdot) =\frac{1}{\eta}K\left(\frac{\cdot}{\eta}\right),
\quad\eta\in (0,1],
\end{displaymath}
and $K\in C^2(\mathbb R;\mathbb R_+)$ is a symmetric $[-1,1]$-supported kernel function. By Menozzi et al. \cite{MPZ21}, Theorem 1.2, the probability distribution of $X_t$ ($t\in (0,T]$) has a continuously differentiable density $f_t$ with respect to Lebesgue's measure on $\mathbb R$ such that $(t,x)\mapsto f_t'(x)$ belongs to $\mathbb L^1([t_0,T]\times\mathbb R)$ when $t_0 > 0$, but not when $t_0 = 0$. So, in order to establish a risk bound on $\widehat f_{\eta}$ (see Lemma \ref{risk_bound_NW_estimator_a}), which is a nonparametric estimator of
\begin{displaymath}
f(\cdot) =\frac{1}{T_0}\int_{t_0}^{T}f_s(\cdot)ds,
\end{displaymath}
the condition $t_0\in (0,T)$ is required. In order to take into account both the regularity and the monotonicity conditions on $a$, our paper deals with the following continuously differentiable and strictly decreasing estimator of $a_{|I_0}$:
\begin{displaymath}
\widehat a_{\ell,h,\eta}(\cdot) =
a(\texttt r_0 +\varepsilon) +
\int_{a(\texttt r_0 +\varepsilon)}^{a(\texttt l_0 -\varepsilon)}\int_{-\infty}^{-\cdot}
K_{\ell}(-\widehat{\mathfrak a}_{h,\eta}^{-1}(z) - y)dydz,
\end{displaymath}
where $I_0 = [\texttt l_0,\texttt r_0]$ with $\texttt l_0,\texttt r_0\in\mathbb R$ satisfying $\texttt l_0 <\texttt r_0$,
\begin{displaymath}
\widehat{\mathfrak a}_{h,\eta}^{-1}(\cdot) =
\texttt l_0 - 2\varepsilon +
\int_{\texttt l_0 - 2\varepsilon}^{\texttt r_0 + 2\varepsilon}\int_{-\infty}^{-\cdot}
K_h(-\widehat a_{\eta}(z) - y)dydz
\quad {\rm and}\quad
\varepsilon,\ell,h\in (0,1].
\end{displaymath}
Since, in practice, $a(\texttt l_0 -\varepsilon)$ and $a(\texttt r_0 +\varepsilon)$ may be unknown, the estimator $\widetilde a_{\ell,h,\eta}$ - with $a$ replaced by $\widehat a_{\eta}$ in the definition of $\widehat a_{\ell,h,\eta}$ - is also investigated.
\\
\\
In the isotonic nonparametric regression framework with $[0,1]$-supported inputs, the asymptotic normality - with a rate of convergence - has been established for an estimator similar to $\widehat{\mathfrak a}_{h,\eta}$ by H. Dette, N. Neumeyer and K.F. Pilz in \cite{DNP06} (see also Birke and Dette \cite{BD08} and Chen et al. \cite{CGFZ19}). Moreover, for an arbitrary estimator $\widehat\beta$ of a strictly increasing function $\beta\in C^2([0,1];\mathbb R)$, N. Neumeyer shows in \cite{NEUMEYER05} that $\|\widehat\beta_h -\beta\|_{\infty,[0,1]}$ is controlled by $\|\widehat\beta -\beta\|_{\infty,[0,1]}$, where
\begin{displaymath}
\widehat\beta_h(\cdot) =
\beta(0) +
\int_{\beta(0)}^{\beta(1)}\mathbf 1_{\widehat\beta_{h}^{\star}(y)\leqslant\cdot}dy
\quad {\rm and}\quad
\widehat\beta_{h}^{\star}(\cdot) =
\int_{0}^{1}\int_{-\infty}^{\cdot}K_h(\widehat\beta(z) - y)dydz.
\end{displaymath}
In these references, the asymptotic properties of the (strictly) monotone estimator are derived from the existing rate of strong uniform convergence on the underlying estimator, as Mack and Silverman \cite{MS82}, Theorem B in the nonparametric regression framework. Since there is no such result for the Nadaraya-Watson estimator $\widehat a_{\eta}$ of the drift function in Equation (\ref{main_equation}), for an arbitrary integrable estimator $\widehat b$ of a strictly decreasing function $b\in C^1(\mathbb R)$, our Section \ref{section_risk_bounds_universal} deals with an $\mathbb L^1$-risk bound - only depending on the $\mathbb L^1$-risk of $\widehat b$ - on the following estimator of $b_{|I_0}$:
\begin{displaymath}
\widehat b_{\ell,h}(\cdot) =
b(\texttt r_0 +\varepsilon) +
\int_{b(\texttt r_0 +\varepsilon)}^{b(\texttt l_0 -\varepsilon)}\int_{-\infty}^{-\cdot}
K_{\ell}(-\widehat{\mathfrak b}_{h}^{-1}(z) - y)dydz
\end{displaymath}
with
\begin{displaymath}
\widehat{\mathfrak b}_{h}^{-1}(\cdot) =
\texttt l_0 - 2\varepsilon +
\int_{\texttt l_0 - 2\varepsilon}^{\texttt r_0 + 2\varepsilon}\int_{-\infty}^{-\cdot}
K_h(-\widehat b(z) - y)dydz.
\end{displaymath}
Moreover, $\mathbb L^1$-risk bounds are established on the estimator $\widetilde b_{\ell,h}$ - with $b$ replaced by $\widehat b$ in the definition of $\widehat b_{\ell,h}$ - and on the adaptive estimator $\widehat b_{\widehat\ell,\widehat h}$, where
\begin{displaymath}
(\widehat\ell,\widehat h) =
\underset{(\ell,h)\in\mathfrak H_{b}^{2}}{\rm argmin}\left\{
\int_{I_0}|\widehat b_{\ell,h}(x) -\widehat b(x)|dx\right\}
\end{displaymath}
and $\mathfrak H_b$ is an appropriate finite subset of $(0,1]$. So, Section \ref{estimation_drift_main_equation} deals with theoretical guarantees on $\widehat a_{\ell,h,\eta}$ (resp. $\widetilde a_{\ell,h,\eta}$), established thanks to the $\mathbb L^1$-risk bound on $\widehat b_{\ell,h}$ (resp. $\widetilde b_{\ell,h}$) by taking $b = a$ and $\widehat b =\widehat a_{\eta}$. Finally, note that in the same spirit as Neumeyer \cite{NEUMEYER05}, for an arbitrary estimator $\widehat\beta$ of an increasing measurable function $\beta : [0,1]\rightarrow\mathbb R$, Chernozhukov et al. \cite{CFVG09} shows that the integrated $\mathbb L^p$-error ($p\geqslant 1$) of
\begin{displaymath}
\widetilde\beta(\cdot) =
\inf\left\{w\in\mathbb R :
\int_{0}^{1}\mathbf 1_{\widehat\beta(y)\leqslant w}dy\geqslant\cdot\right\}
\quad\textrm{is controlled by that of $\widehat\beta$.}
\end{displaymath}
However, in order to take into account both the continuous differentiability and the strict monotonicity of the drift function $a$ in Equation (\ref{main_equation}), our estimator $\widehat a_{\ell,h,\eta}$ seems more appropriate.
\\
\\
The outline of the paper is as follows: the aforementioned theoretical results on $\widehat b_{\ell,h}$ (resp. $\widehat a_{\ell,h,\eta}$) are established in Section \ref{section_risk_bounds_universal} (resp. Section \ref{estimation_drift_main_equation}), and Section \ref{section_numerical_experiments} shows that $\widehat a_{\ell,h,\eta}$ is also satisfactory on the numerical side.
\\
\\
{\bf Notations and basic definitions:}
\begin{itemize}
 \item For every $\lambda > 0$,
 \begin{displaymath}
 I_{\lambda} = [\texttt l_{\lambda},\texttt r_{\lambda}]
 \quad {\rm with}\quad
 \texttt l_{\lambda} =\texttt l_0 -\lambda
 \quad {\rm and}\quad
 \texttt r_{\lambda} =\texttt r_0 +\lambda.
 \end{displaymath}
 Note that $I_{\lambda}\downarrow I_0 = [\texttt l_0,\texttt r_0]$ when $\lambda\downarrow 0$.
 \item Consider an interval $I\subset\mathbb R$. For every continuous function $\varphi : I\rightarrow\mathbb R$,
 \begin{displaymath}
 \|\varphi\|_{\infty,I} =\sup_{x\in I}|\varphi(x)|.
 \end{displaymath}
 Moreover, the set of continuous (resp. $p$-times continuously differentiable ($p\in\mathbb N^*$)) functions from $I$ into $\mathbb R$ is denoted by $C^0(I;\mathbb R)$ (resp. $C^p(I;\mathbb R)$).
\end{itemize}
%


%
\section{Non-asymptotic $\mathbb L^1$-risk bounds on the universal monotone estimator}\label{section_risk_bounds_universal}
Throughout this section, $b$ is a continuously differentiable function from $I_{2\varepsilon}$ into $\mathbb R$ fulfilling the following assumption.
%


%
\begin{assumption}\label{assumption_b}
There exists $\mathfrak m_b > 0$ such that, for every $x\in I_{2\varepsilon}$, $b'(x)\leqslant -\mathfrak m_b$.
\end{assumption}
\noindent
First, let us establish an (integrated) $\mathbb L^1$-risk bound on the continuously differentiable and strictly decreasing estimator
\begin{displaymath}
\widehat b_{\ell,h}(\cdot) =
b(\texttt r_{\varepsilon}) +
\int_{b(I_{\varepsilon})}\int_{-\infty}^{-\cdot}
K_{\ell}(-\widehat{\mathfrak b}_{h}^{-1}(z) - y)dydz
\quad {\rm of}\quad
b_{|I_0},
\end{displaymath}
where
\begin{displaymath}
\widehat{\mathfrak b}_{h}^{-1}(\cdot) =
\texttt l_{2\varepsilon} +
\int_{I_{2\varepsilon}}\int_{-\infty}^{-\cdot}
K_h(-\widehat b(z) - y)dydz,
\end{displaymath}
and $\widehat b$ is an arbitrary integrable estimator of $b_{|I_{2\varepsilon}}$ whose integrated $\mathbb L^1$-risk is denoted by $\mathcal R(\widehat b)$:
\begin{displaymath}
\mathcal R(\widehat b) =
\int_{I_{2\varepsilon}}\mathbb E(|\widehat b(x) - b(x)|)dx.
\end{displaymath}
%


%
\begin{theorem}\label{risk_bound_monotone_estimator_b}
Under Assumption \ref{assumption_b}, if $\ell,h\in (0,\min\{1,\mathfrak m_b\}\varepsilon)$, then
\begin{displaymath}
\int_{I_0}\mathbb E(|\widehat b_{\ell,h}(x) - b(x)|)dx
\lesssim
\ell + h +\mathcal R(\widehat b).
\end{displaymath}
\end{theorem}
\noindent
Now, let us establish an $\mathbb L^1$-risk bound on the adaptive estimator $\widehat b_{\widehat\ell,\widehat h}$, where
\begin{displaymath}
(\widehat\ell,\widehat h) =
\underset{(\ell,h)\in\mathfrak H_{b}^{2}}{\rm argmin}\left\{
\int_{I_0}|\widehat b_{\ell,h}(x) -\widehat b(x)|dx\right\}
\end{displaymath}
and $\mathfrak H_b$ is a finite subset of $(0,\min\{1,\mathfrak m_b\}\varepsilon)$.
%


%
\begin{proposition}\label{risk_bound_adaptive_monotone_estimator_b}
Under Assumption \ref{assumption_b},
\begin{displaymath}
\int_{I_0}\mathbb E(|\widehat b_{\widehat\ell,\widehat h}(x) - b(x)|)dx
\lesssim
\min_{(\ell,h)\in\mathfrak H_{b}^{2}}\left\{\int_{I_0}\mathbb E(|\widehat b_{\ell,h}(x) - b(x)|)dx\right\} +
\mathcal R(\widehat b).
\end{displaymath}
\end{proposition}
%


%
\begin{remark}\label{remark_risk_bound_monotone_estimator_b}
Let us make a few comments about Theorem \ref{risk_bound_monotone_estimator_b} and Proposition \ref{risk_bound_practical_estimator_b}:
\begin{enumerate}
 \item Assume that the integrated $\mathbb L^1$-risk of $\widehat b =\widehat b_N$ converges to zero when $N\rightarrow\infty$, leading to
 \begin{displaymath}
 \mathcal R(\widehat b_N) <\min\{1,\mathfrak m_b\}\varepsilon
 \quad\textrm{for $N$ large enough.}
 \end{displaymath}
 So, one can take
 \begin{displaymath}
 (\ell,h) = (\ell_N,h_N)
 \quad\textrm{with}\quad
 \ell_N,h_N\asymp
 \mathcal R(\widehat b_N),
 \end{displaymath}
 and then the risk bound on $\widehat b_{\ell,h}$ in Theorem \ref{risk_bound_monotone_estimator_b} is of order $\mathcal R(\widehat b_N)$ as expected.
 \item In the framework of Remark \ref{remark_risk_bound_monotone_estimator_b}.(1), assume that $(\ell_N,h_N)\in\mathfrak H_{b}^{2}$. By Theorem \ref{risk_bound_monotone_estimator_b},
 \begin{displaymath}
 \min_{(\ell,h)\in\mathfrak H_{b}^{2}}\left\{\int_{I_0}\mathbb E(|\widehat b_{\ell,h}(x) - b(x)|)dx\right\}
 \lesssim\mathcal R(\widehat b_N),
 \end{displaymath}
 and then Proposition \ref{risk_bound_adaptive_monotone_estimator_b} leads to
 \begin{displaymath}
 \int_{I_0}\mathbb E(|\widehat b_{\widehat\ell,\widehat h}(x) - b(x)|)dx
 \lesssim\mathcal R(\widehat b_N).
 \end{displaymath}
 In other words, the $\mathbb L^1$-risk of the continuously differentiable and strictly decreasing adaptive estimator $\widehat b_{\widehat\ell,\widehat h}$ remains of same order as that of $\widehat b_N$.
\end{enumerate}
\end{remark}
\noindent
Finally, thanks to Theorem \ref{risk_bound_monotone_estimator_b}, let us establish an $\mathbb L^1$-risk bound on the practical estimator
\begin{displaymath}
\widetilde b_{\ell,h}(\cdot) =
\left(\widehat b(\texttt r_{\varepsilon}) +
\int_{\widehat b(\texttt r_{\varepsilon})}^{\widehat b(\texttt l_{\varepsilon})}\int_{-\infty}^{-\cdot}
K_{\ell}(-\widehat{\mathfrak b}_{h}^{-1}(z) - y)dydz\right)
\mathbf 1_{\Omega_{\widehat b}}
\end{displaymath}
with
\begin{displaymath}
\Omega_{\widehat b} =
\left\{\widehat b(\texttt r_{\varepsilon}) -\widehat b(\texttt l_{\varepsilon})
\leqslant -\frac{\mathfrak m_b}{2}(\texttt r_{\varepsilon} -\texttt l_{\varepsilon})\right\}.
\end{displaymath}
To that purpose, consider
\begin{displaymath}
\mathfrak R(\widehat b) =
\sup_{x\in I_{2\varepsilon}}\mathbb E(|\widehat b(x) - b(x)|).
\end{displaymath}
%


%
\begin{proposition}\label{risk_bound_practical_estimator_b}
Under Assumption \ref{assumption_b}, if $\ell,h\in (0,\min\{1,\mathfrak m_b\}\varepsilon)$, then
\begin{displaymath}
\int_{I_0}\mathbb E(|\widetilde b_{\ell,h}(x) - b(x)|)dx
\lesssim
\ell + h +\mathfrak R(\widehat b).
\end{displaymath}
\end{proposition}
%


%
\begin{remark}\label{remark_risk_bound_practical_estimator_b}
Let us make a few comments about Proposition \ref{risk_bound_practical_estimator_b}:
\begin{enumerate}
 \item Assume that - uniformly on $I_{2\varepsilon}$ - the pointwise $\mathbb L^1$-risk of $\widehat b =\widehat b_N$ converges to zero when $N\rightarrow\infty$, leading to
 \begin{displaymath}
 \mathfrak R(\widehat b_N) <
 \min\{1,\mathfrak m_b\}\varepsilon
 \quad\textrm{for $N$ large enough.}
 \end{displaymath}
 So, one can take
 \begin{displaymath}
 (\ell,h) = (\ell_N,h_N)
 \quad\textrm{with}\quad
 \ell_N,h_N\asymp
 \mathfrak R(\widehat b_N),
 \end{displaymath}
 and then the risk bound on $\widetilde b_{\ell,h}$ in Proposition \ref{risk_bound_practical_estimator_b} is of order $\mathfrak R(\widehat b_N)$.
 \item Under Assumption \ref{assumption_b}, by following the same lines as in the proof of Proposition \ref{risk_bound_adaptive_monotone_estimator_b}, one can establish that
 \begin{equation}\label{remark_risk_bound_monotone_estimator_b_1}
 \int_{I_0}\mathbb E(|\widetilde b_{\widetilde\ell,\widetilde h}(x) - b(x)|)dx
 \lesssim
 \min_{(\ell,h)\in\mathfrak H_{b}^{2}}\left\{\int_{I_0}\mathbb E(|\widetilde b_{\ell,h}(x) - b(x)|)dx\right\} +
 \mathfrak R(\widehat b),
 \end{equation}
 where
 \begin{displaymath}
 (\widetilde\ell,\widetilde h) =
 \underset{(\ell,h)\in\mathfrak H_{b}^{2}}{\rm argmin}\left\{
 \int_{I_0}|\widetilde b_{\ell,h}(x) -\widehat b(x)|dx\right\}.
 \end{displaymath}
 In the framework of Remark \ref{remark_risk_bound_practical_estimator_b}.(1), assume that $(\ell_N,h_N)\in\mathfrak H_{b}^{2}$. By Proposition \ref{risk_bound_practical_estimator_b},
 \begin{displaymath}
 \min_{(\ell,h)\in\mathfrak H_{b}^{2}}\left\{\int_{I_0}\mathbb E(|\widetilde b_{\ell,h}(x) - b(x)|)dx\right\}
 \lesssim\mathfrak R(\widehat b_N),
 \end{displaymath}
 and then Inequality (\ref{remark_risk_bound_monotone_estimator_b_1}) leads to
 \begin{displaymath}
 \int_{I_0}\mathbb E(|\widetilde b_{\widetilde\ell,\widetilde h}(x) - b(x)|)dx
 \lesssim\mathfrak R(\widehat b_N).
 \end{displaymath}
\end{enumerate}
\end{remark}
%


%
\section{Back to the estimation of $a$ in Model (\ref{main_equation})}\label{estimation_drift_main_equation}
First, thanks to Theorem \ref{risk_bound_monotone_estimator_b}, let us establish an $\mathbb L^1$-risk bound on our continuously differentiable and strictly decreasing estimator of $a_{|I_0}$.
%


%
\begin{theorem}\label{risk_bound_monotone_estimator_a}
Under Assumption \ref{assumption_drift_volatility}, if $\ell,h\in (0,\min\{1,\mathfrak m_a\}\varepsilon)$, then
\begin{displaymath}
\int_{I_0}\mathbb E(|\widehat a_{\ell,h,\eta}(x) - a(x)|)dx
\lesssim
\ell + h +\left(\eta^2 +\frac{1}{N\eta}\right)^{\frac{1}{2}}.
\end{displaymath}
\end{theorem}
%


%
\begin{remark}\label{remark_risk_bound_monotone_estimator_a}
Let us make a few comments about Theorem \ref{risk_bound_monotone_estimator_a}:
\begin{enumerate}
 \item For $N$ large enough, if $\eta,\ell,h\asymp N^{-\frac{1}{3}}$, then the $\mathbb L^1$-rate of $\widehat a_{\ell,h,\eta}$ is of order $N^{-\frac{1}{3}}$ as expected.
 \item In practice, a suitable way to select $\eta$ from data is the extension of the leave-one-out cross-validation (LooCV) procedure provided in Marie and Rosier \cite{MR23}, Section 5.2:
 \begin{equation}\label{remark_risk_bound_monotone_estimator_a_1}
 \widehat\eta =\underset{\eta\in\mathfrak H}{\rm argmin}\left\{
 \sum_{i = 1}^{N}\left(
 \int_{t_0}^{T}\widehat a_{\eta,i}(X_{s}^{i})^2ds -
 2\int_{t_0}^{T}\widehat a_{\eta,i}(X_{s}^{i})dX_{s}^{i}
 \right)\right\},
 \end{equation}
 where $\mathfrak H$ is a finite subset of $(0,1]$ and, for every $i\in\{1,\dots,N\}$ and $\eta\in (0,1]$,
 \begin{displaymath}
 \widehat a_{\eta,i}(\cdot) =
 \frac{\widehat{af}_{\eta,i}(\cdot)}{\widehat f_{\eta}(\cdot)}
 \quad\textrm{with}\quad
 \widehat{af}_{\eta,i}(\cdot) =
 \frac{1}{NT_0}\sum_{k\neq i}\int_{t_0}^{T}K_{\eta}(X_{s}^{k} -\cdot)dX_{s}^{k}.
 \end{displaymath}
 \item Consider
 \begin{equation}\label{remark_risk_bound_monotone_estimator_a_2}
 \widehat{\bf w} =
 \underset{\mathbf w\in\mathfrak H_{a}^{2}}{\rm argmin}\left\{
 \int_{I_0}|\widehat a_{\mathbf w,\widehat\eta}(x) -\widehat a_{\widehat\eta}(x)|dx\right\},
 \end{equation}
 where $\widehat a_{\widehat\eta}$ is an adaptive version of the copies-based Nadaraya-Watson estimator of $a$. By Theorem \ref{risk_bound_monotone_estimator_b} and Proposition \ref{risk_bound_adaptive_monotone_estimator_b} with $b = a$ and $\widehat b =\widehat a_{\widehat\eta}$,
 \begin{equation}\label{remark_risk_bound_monotone_estimator_a_3}
 \int_{I_0}\mathbb E(|\widehat a_{\widehat{\bf w},\widehat\eta}(x) - a(x)|)dx
 \lesssim
 \min_{(\ell,h)\in\mathfrak H_{a}^{2}}\{
 \ell + h +\mathcal R(\widehat a_{\widehat\eta})\} +
 \mathcal R(\widehat a_{\widehat\eta}).
 \end{equation}
 So, even if there is no $\mathbb L^1$-risk bound on $\widehat a_{\widehat\eta}$ - as when $\widehat\eta$ is selected from data via the LooCV procedure presented above - Inequality (\ref{remark_risk_bound_monotone_estimator_a_3}) says that if $\mathcal R(\widehat a_{\widehat\eta})$ belongs to $\mathfrak H_a$, then the $\mathbb L^1$-risk of $\widehat a_{\widehat{\bf w},\widehat\eta}$ is controlled by $\mathcal R(\widehat a_{\widehat\eta})$. Moreover, note that by replacing $\widehat a_{\widehat\eta}$ by the 2-bandwidths PCO-adaptive Nadaraya-Watson estimator investigated in Marie and Rosier \cite{MR23}, Section 5.1, one can establish an oracle inequality on the corresponding version of our continuously differentiable and strictly decreasing estimator of $a_{|I_0}$ thanks to Proposition \ref{risk_bound_adaptive_monotone_estimator_b} and Marie and Rosier \cite{MR23}, Corollary 2. This is left to the reader because the experiments presented in Section \ref{section_numerical_experiments} are based on the $1$-bandwidth LooCV-adaptive Nadaraya-Watson estimator of $a$, which is more satisfactory on the numerical side than the 2-bandwidths PCO-adaptive one.
\end{enumerate}
\end{remark}
\noindent
Now, since $a(\texttt l_0 -\varepsilon)$ and $a(\texttt r_0 +\varepsilon)$ may be unknown in practice, let us establish an $\mathbb L^1$-risk bound on the practical estimator
\begin{displaymath}
\widetilde a_{\ell,h,\eta}(\cdot) =
\left(\widehat a_{\eta}(\texttt r_{\varepsilon}) +
\int_{\widehat a_{\eta}(\texttt r_{\varepsilon})}^{\widehat a_{\eta}(\texttt l_{\varepsilon})}
\int_{-\infty}^{-\cdot}K_{\ell}(-\widehat{\mathfrak a}_{h,\eta}^{-1}(z) - y)dydz\right)
\mathbf 1_{\Omega_{\widehat a_{\eta}}}
\end{displaymath}
with
\begin{displaymath}
\Omega_{\widehat a_{\eta}} =
\left\{\widehat a_{\eta}(\texttt r_{\varepsilon}) -\widehat a_{\eta}(\texttt l_{\varepsilon})
\leqslant -\frac{\mathfrak m_a}{2}(\texttt r_{\varepsilon} -\texttt l_{\varepsilon})\right\}.
\end{displaymath}
%


%
\begin{proposition}\label{risk_bound_practical_estimator_a}
Under Assumption \ref{assumption_drift_volatility}, if $\ell,h\in (0,\min\{1,\mathfrak m_a\}\varepsilon)$, then
\begin{displaymath}
\int_{I_0}\mathbb E(|\widetilde a_{\ell,h,\eta}(x) - a(x)|)dx
\lesssim
\ell + h +
\left(\eta^2 +\frac{1}{N\eta}\right)^{\frac{1}{2}}.
\end{displaymath}
\end{proposition}
%


%
\begin{remark}\label{remark_risk_bound_practical_estimator_a}
Let us make a few comments about Proposition \ref{risk_bound_practical_estimator_a}:
\begin{enumerate}
 \item For $N$ large enough, if $\eta,\ell,h\asymp N^{-\frac{1}{3}}$, then the $\mathbb L^1$-rate of $\widetilde a_{\ell,h,\eta}$ is of order $N^{-\frac{1}{3}}$ as that of $\widehat a_{\ell,h,\eta}$.
 \item Consider
 \begin{displaymath}
 \widetilde{\bf w} =
 \underset{\mathbf w\in\mathfrak H_{a}^{2}}{\rm argmin}\left\{
 \int_{I_0}|\widetilde a_{\mathbf w,\widehat\eta}(x) -\widehat a_{\widehat\eta}(x)|dx\right\},
 \end{displaymath}
 where $\widehat a_{\widehat\eta}$ is an adaptive version of the copies-based Nadaraya-Watson estimator of $a$. By Proposition \ref{risk_bound_practical_estimator_b} and Inequality (\ref{remark_risk_bound_monotone_estimator_b_1}) with $b = a$ and $\widehat b =\widehat a_{\widehat\eta}$,
 \begin{equation}\label{remark_risk_bound_practical_estimator_a_1}
 \int_{I_0}\mathbb E(|\widetilde a_{\widetilde{\bf w},\widehat\eta}(x) - a(x)|)dx
 \lesssim
 \min_{(\ell,h)\in\mathfrak H_{a}^{2}}\{
 \ell + h +\mathfrak R(\widehat a_{\widehat\eta})\} +
 \mathfrak R(\widehat a_{\widehat\eta}).
 \end{equation}
 So, even if there is no $\mathbb L^1$-risk bound on $\widehat a_{\widehat\eta}$ - as when $\widehat\eta$ is selected from data via the LooCV procedure presented in Remark \ref{remark_risk_bound_monotone_estimator_a}.(2) - Inequality (\ref{remark_risk_bound_practical_estimator_a_1}) says that if $\mathfrak R(\widehat a_{\widehat\eta})$ belongs to $\mathfrak H_a$, then the $\mathbb L^1$-risk of $\widetilde a_{\widetilde{\bf w},\widehat\eta}$ is controlled by $\mathfrak R(\widehat a_{\widehat\eta})$.
\end{enumerate}
\end{remark}
%


%
\section{Numerical experiments}\label{section_numerical_experiments}
\begin{figure}[!h]
\centering
\includegraphics[scale=0.55]{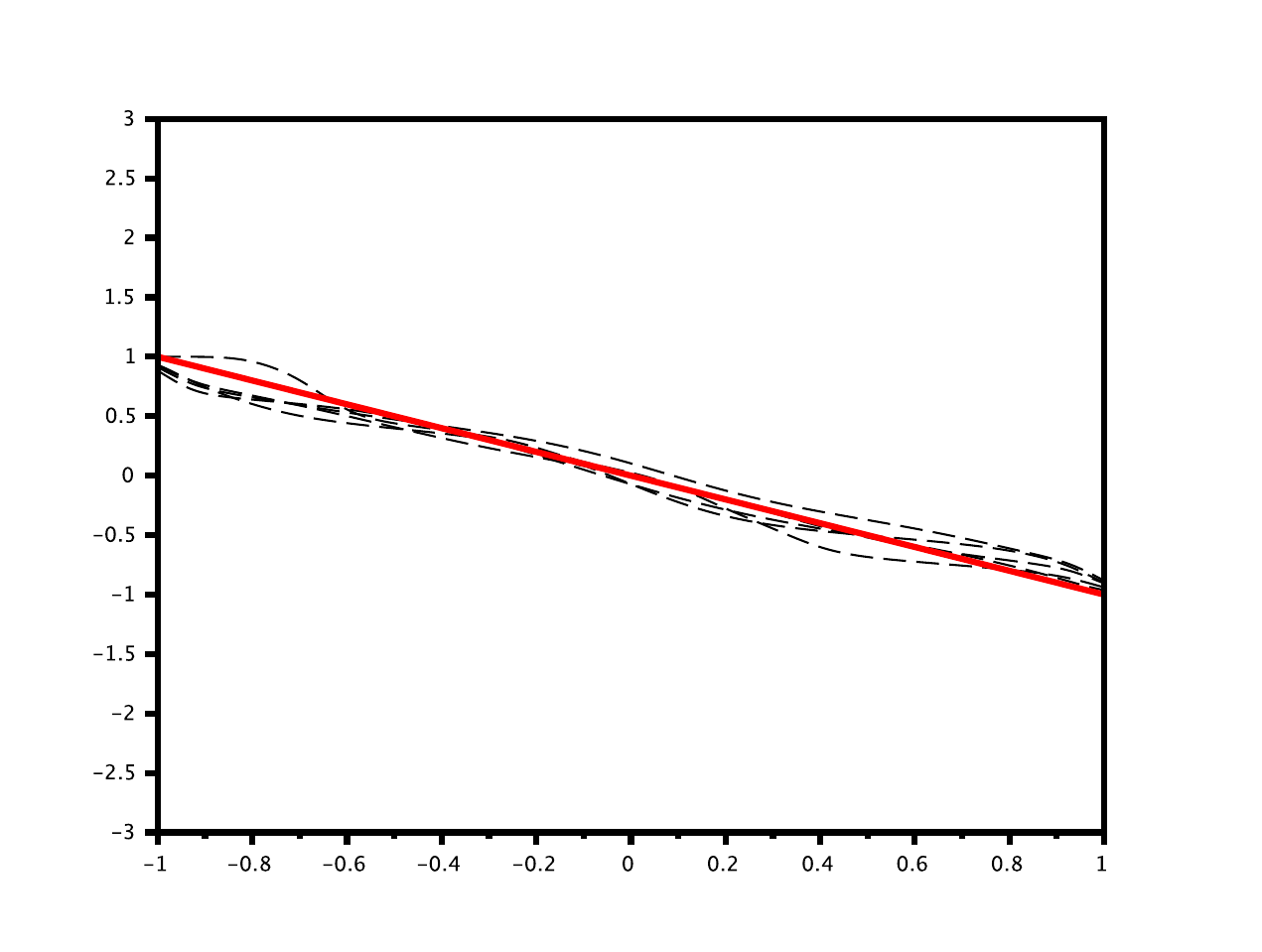}
\caption{Plots of 5 adaptive strictly decreasing estimations (black dashed lines) of $a$ (red line) for Model (A).}
\label{adaptive_estimations_Model_A}
\end{figure}

\begin{figure}[!h]
\centering
\includegraphics[scale=0.55]{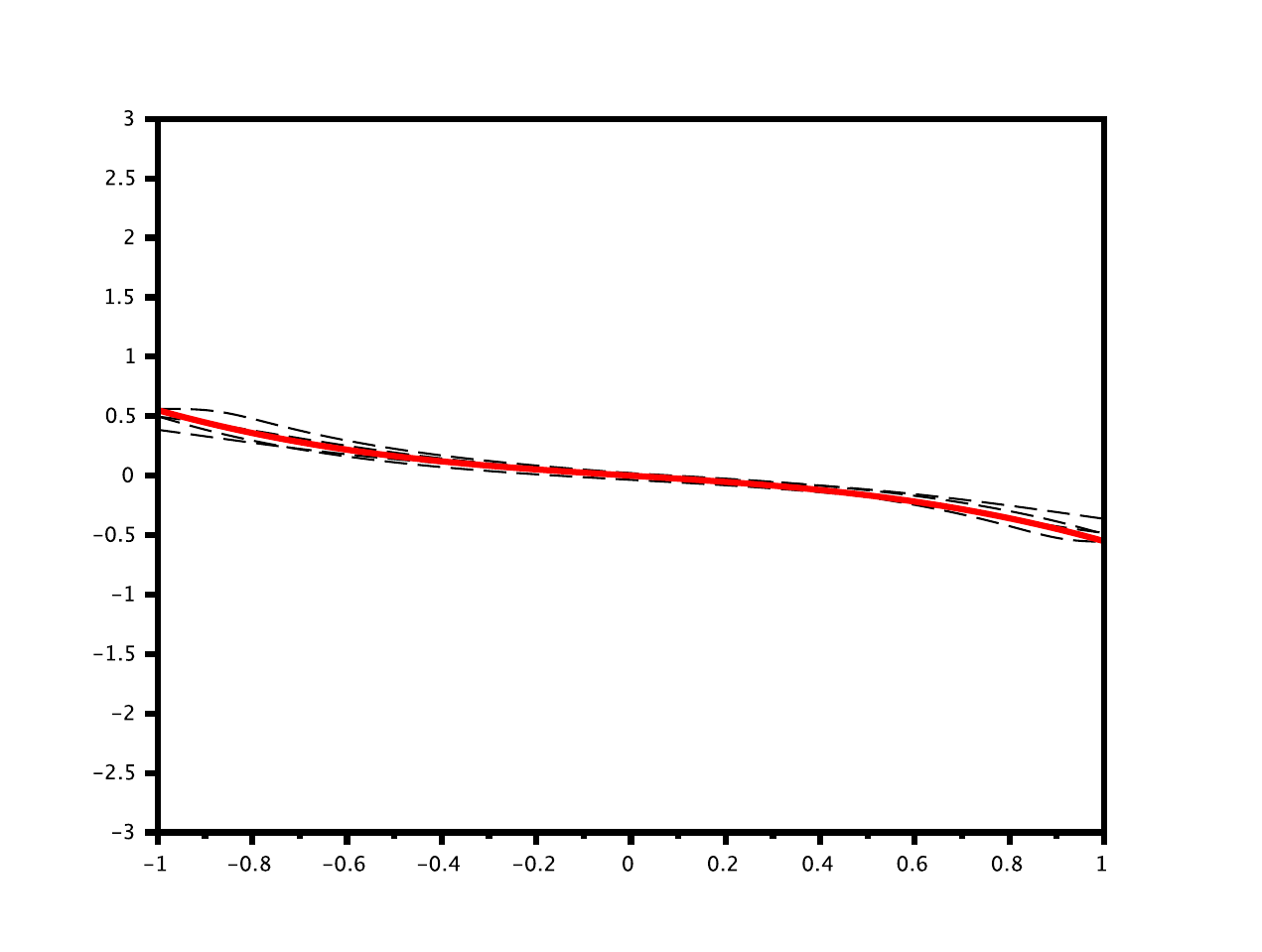}
\caption{Plots of 5 adaptive strictly decreasing estimations (black dashed lines) of $a$ (red line) for Model (B).}
\label{adaptive_estimations_Model_B}
\end{figure}

In this section, numerical experiments on the adaptive strictly decreasing estimator $\widehat a_{\widehat{\bf w},\widehat\eta}$, where $\widehat\eta$ (resp. $\widehat{\bf w}$) is selected from data thanks to the LooCV criterion (\ref{remark_risk_bound_monotone_estimator_a_1}) (resp. to (\ref{remark_risk_bound_monotone_estimator_a_2})), are presented for the two following models:
\begin{itemize}
 \item[(A)] $\displaystyle{X_t = 0.5 -\int_{0}^{t}X_sds + W_t}$, and
 \item[(B)] $\displaystyle{X_t = 0.5 +\int_{0}^{t}\left(\sin\left(\frac{5X_s}{4}\right) -\frac{3X_s}{2}\right)ds + W_t}$.
\end{itemize}
For each model, $\widehat a_{\widehat{\bf w},\widehat\eta}$ is computed on $I_0 = [-1,1]$ from $N = 100$ paths of the process $X$ observed along the dissection $\{\lambda T/n\textrm{ $;$ }\lambda = 0,\dots,n\}$ of $[0,T]$, where $n = 50$, $T = 5$, $K$ is the standard normal density function, $\varepsilon = 0.01$, $\widehat\eta$ is selected in $\mathfrak H =\{0.05\lambda\textrm{ $;$ }\lambda = 1,\dots,35\}$, and $\widehat{\bf w}$ is selected in $\mathfrak H^2$. This experiment is repeated 100 times, and 5 adaptive strictly decreasing estimations of the linear (resp. nonlinear) function $a$ are plotted on Figure \ref{adaptive_estimations_Model_A} (resp. Figure \ref{adaptive_estimations_Model_B}) for Model (A) (resp. Model (B)).
\newline
For both Models (A) and (B), the means and standard deviations of the integrated errors of $\widehat a_{\widehat{\bf w},\widehat\eta}$ and $\widehat a_{\widehat\eta}$ are stored in Table \ref{table_L1_risks}. These are all small, which is visible on Figures \ref{adaptive_estimations_Model_A} and \ref{adaptive_estimations_Model_B}. Moreover, the mean of the integrated error of $\widehat a_{\widehat{\bf w},\widehat\eta}$ is of same order as (resp. lower than) that of $\widehat a_{\widehat\eta}$ for Model (A) (resp. Model (B)), while the standard deviations of the integrated errors of $\widehat a_{\widehat{\bf w},\widehat\eta}$ and $\widehat a_{\widehat\eta}$ are of same order for both Models (A) and (B).
\begin{table}[h!]
\begin{center}
\begin{tabular}{cll}
\hline
 Model & Strict. decreasing estim. & Nadaraya-Watson estim.\\
\hline
 (A) & $0.162_{(0.046)}$ & $0.197_{(0.043)}$\\
 (B) & $0.070_{(0.030)}$ & $0.162_{(0.043)}$\\
\hline
\end{tabular}
\medskip
\caption{Mean integrated errors of $\widehat a_{\widehat{\bf w},\widehat\eta}$ and $\widehat a_{\widehat\eta}$ (with standard deviations) computed over 100 repetitions for Models (A) and (B).}\label{table_L1_risks}
\end{center}
\end{table}
%


%
\section{Proofs}\label{section_proofs}
%


%
\subsection{Proof of Theorem \ref{risk_bound_monotone_estimator_b}}\label{section_proof_risk_bound_monotone_estimator_b}
Consider
\begin{displaymath}
\mathfrak b_{h}^{-1}(\cdot) =
\texttt l_{2\varepsilon} +
\int_{I_{2\varepsilon}}\int_{-\infty}^{-\cdot}K_h(-b(z) - y)dydz,
\end{displaymath}
which is a strictly decreasing approximation of $b^{-1}$ on $b(I_{\varepsilon})$, and
\begin{displaymath}
b_{\ell}(\cdot) =
b(\texttt r_{\varepsilon}) +
\int_{b(I_{\varepsilon})}\int_{-\infty}^{-\cdot}K_{\ell}(-b^{-1}(z) - y)dydz,
\end{displaymath}
which is a strictly decreasing approximation of $b$ on $I_0$. Indeed, for every $w\in b(I_{\varepsilon})$,
\begin{eqnarray*}
 \mathfrak b_{h}^{-1}(w) & = &
 \texttt l_{2\varepsilon} +
 \int_{I_{2\varepsilon}}
 \int_{-\infty}^{\infty}K_h(-b(z) - y)\mathbf 1_{y\leqslant -w}dydz\\
 & = &
 \texttt l_{2\varepsilon} +
 \int_{I_{2\varepsilon}}(K_h\star\mathbf 1_{(-\infty,-w]})(-b(z))dz
 \xrightarrow[h\rightarrow 0^+]{}
 \texttt l_{2\varepsilon} +
 \int_{I_{2\varepsilon}}
 \mathbf 1_{z\leqslant b^{-1}(w)}dz = b^{-1}(w),
\end{eqnarray*}
and for every $x\in I_0$,
\begin{eqnarray*}
 b_{\ell}(x) & = &
 b(\texttt r_{\varepsilon}) +
 \int_{b(I_{\varepsilon})}
 \int_{-\infty}^{\infty}K_{\ell}(-b^{-1}(z) - y)\mathbf 1_{y\leqslant -x}dydz\\
 & = &
 b(\texttt r_{\varepsilon}) +
 \int_{b(I_{\varepsilon})}(K_{\ell}\star\mathbf 1_{(-\infty,-x]})(-b^{-1}(z))dz
 \xrightarrow[\ell\rightarrow 0^+]{}
 b(\texttt r_{\varepsilon}) +
 \int_{b(I_{\varepsilon})}
 \mathbf 1_{z\leqslant b(x)}dz = b(x).
\end{eqnarray*}
The proof of Theorem \ref{risk_bound_monotone_estimator_b} relies on the following lemma.
%


%
\begin{lemma}\label{bounds_approximations}
Under Assumption \ref{assumption_b},
\begin{enumerate}
 \item If $h <\min\{b({\tt l}_{2\varepsilon}) - b({\tt l}_{\varepsilon}),b({\tt r}_{\varepsilon}) - b({\tt r}_{2\varepsilon})\}$, then $\displaystyle{
 \sup_{w\in b(I_{\varepsilon})}|\mathfrak b_{h}^{-1}(w) - b^{-1}(w)|\lesssim h}$.
 \item If $\ell <\varepsilon$, then $\displaystyle{\sup_{x\in I_0}|b_{\ell}(x) - b(x)|\lesssim\ell}$.
\end{enumerate}
\end{lemma}
\noindent
The proof of Lemma \ref{bounds_approximations} is postponed to Section \ref{section_proof_bounds_approximations}. First, consider
\begin{displaymath}
\Delta_h(w) :=
\int_{I_{2\varepsilon}}\int_{-\infty}^{-w}
(K_h(-\widehat b(z) - y) - K_h(-b(z) - y))dydz,
\quad w\in\mathbb R,
\end{displaymath}
and note that
\begin{displaymath}
\widehat{\mathfrak b}_{h}^{-1}(w) -\mathfrak b_{h}^{-1}(w) =\Delta_h(w),
\quad\forall w\in b(I_{\varepsilon}).
\end{displaymath}
By Taylor's formula with integral remainder, for every $w\in\mathbb R$,
\begin{eqnarray*}
 \Delta_h(w)
 & = &
 \int_{I_{2\varepsilon}}\int_{-\infty}^{-w}
 (b(z) -\widehat b(z))\int_{0}^{1}K_h'(-b(z) +\xi (b(z) -\widehat b(z)) - y)d\xi dydz\\
 & = &
 -\int_{I_{2\varepsilon}}(b(z) -\widehat b(z))
 \int_{0}^{1}K_h(-b(z) +\xi (b(z) -\widehat b(z)) + w)d\xi dz.
\end{eqnarray*}
So,
\begin{eqnarray*}
 \int_{-\infty}^{\infty}\mathbb E(|\Delta_h(w)|)dw
 & \leqslant &
 \int_{I_{2\varepsilon}}
 \mathbb E\left(|\widehat b(z) - b(z)|
 \int_{0}^{1}\int_{-\infty}^{\infty}|K_h(-b(z) +\xi (b(z) -\widehat b(z)) + w)|dwd\xi\right)dz\\
 & \leqslant &
 \|K\|_1
 \int_{I_{2\varepsilon}}\mathbb E(|\widehat b(z) - b(z)|)dz =\mathcal R(\widehat b).
\end{eqnarray*}
Moreover, by the mean value theorem, and by Assumption \ref{assumption_b},
\begin{displaymath}
b({\tt l}_{2\varepsilon}) - b({\tt l}_{\varepsilon})\geqslant
\mathfrak m_b(\texttt l_{\varepsilon} -\texttt l_{2\varepsilon}) =
\mathfrak m_b\varepsilon
\quad {\rm and}\quad
b({\tt r}_{\varepsilon}) - b({\tt r}_{2\varepsilon})\geqslant
\mathfrak m_b(\texttt r_{2\varepsilon} -\texttt r_{\varepsilon}) =
\mathfrak m_b\varepsilon,
\end{displaymath}
leading to
\begin{displaymath}
h <\mathfrak m_b\varepsilon\leqslant
\min\{b({\tt l}_{2\varepsilon}) - b({\tt l}_{\varepsilon}),
b({\tt r}_{\varepsilon}) - b({\tt r}_{2\varepsilon})\}.
\end{displaymath}
Thus, by Lemma \ref{bounds_approximations}.(1),
\begin{eqnarray}
 \int_{b(I_{\varepsilon})}\mathbb E(|\widehat{\mathfrak b}_{h}^{-1}(w) - b^{-1}(w)|)dw
 & \leqslant &
 \int_{b(I_{\varepsilon})}|\mathfrak b_{h}^{-1}(w) - b^{-1}(w)|dw +
 \int_{b(I_{\varepsilon})}\mathbb E(|\Delta_h(w)|)dw
 \nonumber\\
 \label{risk_bound_monotone_estimator_b_1}
 & \lesssim &
 h +\mathcal R(\widehat b).
\end{eqnarray}
Now, consider
\begin{displaymath}
\delta_{\ell,h}(x) :=
\int_{b(I_{\varepsilon})}\int_{-\infty}^{-x}
(K_{\ell}(-\widehat{\mathfrak b}_{h}^{-1}(z) - y) - K_{\ell}(-b^{-1}(z) - y))dydz,
\quad x\in\mathbb R,
\end{displaymath}
and note that
\begin{displaymath}
\widehat b_{\ell,h}(x) - b_{\ell}(x) =\delta_{\ell,h}(x),
\quad\forall x\in I_0.
\end{displaymath}
By Taylor's formula with integral remainder, for every $x\in\mathbb R$,
\begin{eqnarray*}
 \delta_{\ell,h}(x)
 & = &
 \int_{b(I_{\varepsilon})}\int_{-\infty}^{-x}
 (b^{-1}(z) -\widehat{\mathfrak b}_{h}^{-1}(z))\int_{0}^{1}
 K_{\ell}'(-b^{-1}(z) +\xi (b^{-1}(z) -\widehat{\mathfrak b}_{h}^{-1}(z)) - y)d\xi dydz\\
 & = &
 -\int_{b(I_{\varepsilon})}(b^{-1}(z) -\widehat{\mathfrak b}_{h}^{-1}(z))
 \int_{0}^{1}K_{\ell}(-b^{-1}(z) +\xi (b^{-1}(z) -\widehat{\mathfrak b}_{h}^{-1}(z)) + x)d\xi dz.
\end{eqnarray*}
So, by Inequality (\ref{risk_bound_monotone_estimator_b_1}),
\begin{eqnarray*}
 \int_{-\infty}^{\infty}\mathbb E(|\delta_{\ell,h}(x)|)dx
 & \leqslant &
 \int_{b(I_{\varepsilon})}
 \mathbb E(|\widehat{\mathfrak b}_{h}^{-1}(z) - b^{-1}(z)|\\
 & &
 \hspace{2cm}\times
 \int_{0}^{1}\int_{-\infty}^{\infty}
 |K_{\ell}(-b^{-1}(z) +\xi (b^{-1}(z) -\widehat{\mathfrak b}_{h}^{-1}(z)) + x)|dxd\xi)dz\\
 & = &
 \int_{b(I_{\varepsilon})}\mathbb E(|\widehat{\mathfrak b}_{h}^{-1}(z) - b^{-1}(z)|)dz
 \lesssim
 h +\mathcal R(\widehat b).
\end{eqnarray*}
Therefore, by Lemma \ref{bounds_approximations}.(2),
\begin{eqnarray*}
 \int_{I_0}\mathbb E(|\widehat b_{\ell,h}(x) - b(x)|)dx
 & \leqslant &
 \int_{I_0}|b_{\ell}(x) - b(x)|dx +
 \int_{I_0}\mathbb E(|\delta_{\ell,h}(x)|)dx\\
 & \lesssim &
 \ell + h +\mathcal R(\widehat b).
\end{eqnarray*}
%


%
\subsubsection{Proof of Lemma \ref{bounds_approximations}:}\label{section_proof_bounds_approximations}
\begin{enumerate}
 \item Consider $w\in b(I_{\varepsilon})$, and note that for every $y\in [-1,1]$,
 \begin{displaymath}
 \texttt l_{2\varepsilon}\leqslant
 b^{-1}(b(\texttt l_{\varepsilon}) + h)\leqslant
 b^{-1}(w - hy)\leqslant
 b^{-1}(b(\texttt r_{\varepsilon}) - h)\leqslant
 \texttt r_{2\varepsilon}.
 \end{displaymath}
 Obviously, $b^{-1}(w)$ also belongs to $I_{2\varepsilon}$. Then,
 \begin{eqnarray*}
  |\mathfrak b_{h}^{-1}(w) - b^{-1}(w)| & = &
  \left|\int_{I_{2\varepsilon}}
  \int_{-\infty}^{\infty}K_h(-b(z) - y)\mathbf 1_{y\leqslant -w}dydz -
  \int_{I_{2\varepsilon}}\mathbf 1_{-b(z)\leqslant -w}dz\right|\\
  & = &
  \left|\int_{I_{2\varepsilon}}
  \int_{-\infty}^{\infty}K(y)(\mathbf 1_{-b(z) - hy\leqslant -w} -\mathbf 1_{-b(z)\leqslant -w})dydz\right|\\
  & = &
  \left|\int_{I_{2\varepsilon}}
  \int_{-\infty}^{\infty}K(y)(\mathbf 1_{z\leqslant b^{-1}(w - hy)} -\mathbf 1_{z\leqslant b^{-1}(w)})dydz\right|\\
  & = &
  \left|\int_{-\infty}^{\infty}K(y)(b^{-1}(w - hy) - b^{-1}(w))dy\right|
  \leqslant
  h\|(b^{-1})'\|_{\infty,b(I_{2\varepsilon})}\int_{-\infty}^{\infty}|y|K(y)dy.
 \end{eqnarray*}
 \item Consider $x\in I_0$, and note that for every $y\in [-1,1]$,
 \begin{displaymath}
 b(\texttt r_{\varepsilon})\leqslant
 b(\texttt r_0 +\ell)\leqslant
 b(x -\ell y)\leqslant
 b(\texttt l_0 -\ell)\leqslant
 b(\texttt l_{\varepsilon}).
 \end{displaymath}
 Obviously, $b(x)$ also belongs to $b(I_{\varepsilon})$. Then,
 \begin{eqnarray*}
  |b_{\ell}(x) - b(x)| & = &
  \left|\int_{b(I_{\varepsilon})}
  \int_{-\infty}^{\infty}K_{\ell}(-b^{-1}(z) - y)\mathbf 1_{y\leqslant -x}dydz -
  \int_{b(I_{\varepsilon})}\mathbf 1_{-b^{-1}(z)\leqslant -x}dz\right|\\
  & = &
  \left|\int_{b(I_{\varepsilon})}
  \int_{-\infty}^{\infty}K(y)(\mathbf 1_{z\leqslant b(x -\ell y)} -\mathbf 1_{z\leqslant b(x)})dydz\right|
  \leqslant
  \ell\|b'\|_{\infty,I_{\varepsilon}}\int_{-\infty}^{\infty}|y|K(y)dy.
 \end{eqnarray*}
\end{enumerate}
%


%
\subsection{Proof of Proposition \ref{risk_bound_adaptive_monotone_estimator_b}}\label{section_proof_risk_bound_adaptive_monotone_estimator_b}
By the definition of $(\widehat\ell,\widehat h)$,
\begin{eqnarray*}
 \int_{I_0}|\widehat b_{\widehat\ell,\widehat h}(x) -\widehat b(x)|dx
 & = &
 \min_{(\ell,h)\in\mathfrak H_b}\left\{\int_{I_0}|\widehat b_{\ell,h}(x) -\widehat b(x)|dx\right\}\\
 & \leqslant &
 \min_{(\ell,h)\in\mathfrak H_b}\left\{\int_{I_0}|\widehat b_{\ell,h}(x) - b(x)|dx +
 \int_{I_0}|\widehat b(x) - b(x)|dx\right\}\\
 & = &
 \min_{(\ell,h)\in\mathfrak H_b}\left\{\int_{I_0}|\widehat b_{\ell,h}(x) - b(x)|dx\right\} +
 \int_{I_0}|\widehat b(x) - b(x)|dx,
\end{eqnarray*}
leading to
\begin{displaymath}
\int_{I_0}\mathbb E(|\widehat b_{\widehat\ell,\widehat h}(x) -\widehat b(x)|)dx
\leqslant
\min_{(\ell,h)\in\mathfrak H_b}\left\{\int_{I_0}\mathbb E(|\widehat b_{\ell,h}(x) - b(x)|)dx\right\} +\mathcal R(\widehat b).
\end{displaymath}
Thus,
\begin{eqnarray*}
 \int_{I_0}\mathbb E(|\widehat b_{\widehat\ell,\widehat h}(x) - b(x)|)dx
 & \leqslant &
 \int_{I_0}\mathbb E(|\widehat b_{\widehat\ell,\widehat h}(x) -\widehat b(x)|)dx +\mathcal R(\widehat b)\\
 & \leqslant &
 \min_{(\ell,h)\in\mathfrak H_b}\left\{\int_{I_0}\mathbb E(|\widehat b_{\ell,h}(x) - b(x)|)dx\right\} + 2\mathcal R(\widehat b).
\end{eqnarray*}
%


%
\subsection{Proof of Proposition \ref{risk_bound_practical_estimator_b}}\label{section_proof_risk_bound_practical_estimator_b}
First, consider
\begin{displaymath}
\widehat b_{\ell,h}^{\star}(\cdot) =
\widehat b(\texttt r_{\varepsilon}) +
\int_{\widehat b(\texttt r_{\varepsilon})}^{\widehat b(\texttt l_{\varepsilon})}\int_{-\infty}^{-\cdot}
K_{\ell}(-\widehat{\mathfrak b}_{h}^{-1}(z) - y)dydz,
\end{displaymath}
and note that
\begin{eqnarray}
 \int_{I_0}\mathbb E(
 |\widetilde b_{\ell,h}(x) - b(x)|)dx & = &
 \mathbb E\left(\mathbf 1_{\Omega_{\widehat b}}
 \int_{I_0}|\widehat b_{\ell,h}^{\star}(x) - b(x)|dx\right) +
 \mathbb E\left(\mathbf 1_{\Omega_{\widehat b}^{c}}
 \int_{I_0}|b(x)|dx\right)
 \nonumber\\
 \label{risk_bound_practical_estimator_b_1}
 & \leqslant &
 \int_{I_0}\mathbb E(|\widehat b_{\ell,h}(x) - b(x)|)dx +
 (\texttt r_0 -\texttt l_0)\|b\|_{\infty,I_0}\mathbb P(\Omega_{\widehat b}^{c}) + R_{\varepsilon}
\end{eqnarray}
with
\begin{displaymath}
R_{\varepsilon} =
\mathbb E\left(\mathbf 1_{\Omega_{\widehat b}}
\int_{I_0}|\widehat b_{\ell,h}^{\star}(x) -\widehat b_{\ell,h}(x)|dx\right).
\end{displaymath}
Now, let us establish suitable bounds on $R_{\varepsilon}$ and $\mathbb P(\Omega_{\widehat b}^{c})$.
\\
\\
{\bf Bound on $R_{\varepsilon}$.} Consider the $[0,1]$-valued function
\begin{displaymath}
\mathcal K : w\in\mathbb R\longmapsto\int_{-\infty}^{w}K(y)dy.
\end{displaymath}
For every $x\in I_0$,
\begin{eqnarray*}
 |\widehat b_{\ell,h}^{\star}(x) -\widehat b_{\ell,h}(x)|\mathbf 1_{\Omega_{\widehat b}}
 & = &
 \left|\widehat b(\texttt l_{\varepsilon}) -
 \int_{\widehat b(\texttt r_{\varepsilon})}^{\widehat b(\texttt l_{\varepsilon})}
 \mathcal K\left(\frac{-\widehat{\mathfrak b}_{h}^{-1}(z) + x}{\ell}\right)dz\right.\\
 & &
 \hspace{2cm}\left. -
 b(\texttt l_{\varepsilon}) +
 \int_{b(\texttt r_{\varepsilon})}^{b(\texttt l_{\varepsilon})}
 \mathcal K\left(\frac{-\widehat{\mathfrak b}_{h}^{-1}(z) + x}{\ell}\right)dz\right|
 \mathbf 1_{\Omega_{\widehat b}}\\
 & \leqslant &
 |\widehat b(\texttt l_{\varepsilon}) - b(\texttt l_{\varepsilon})| +
 \left|\int_{\widehat b(\texttt r_{\varepsilon})}^{b(\texttt r_{\varepsilon})}
 \mathcal K\left(\frac{-\widehat{\mathfrak b}_{h}^{-1}(z) + x}{\ell}\right)dz\right| +
 \left|\int_{b(\texttt l_{\varepsilon})}^{\widehat b(\texttt l_{\varepsilon})}
 \mathcal K\left(\frac{-\widehat{\mathfrak b}_{h}^{-1}(z) + x}{\ell}\right)dz\right|\\
 & \leqslant &
 2|\widehat b(\texttt l_{\varepsilon}) - b(\texttt l_{\varepsilon})| +
 |\widehat b(\texttt r_{\varepsilon}) - b(\texttt r_{\varepsilon})|.
\end{eqnarray*}
Thus,
\begin{displaymath}
R_{\varepsilon}\leqslant
3(\texttt r_0 -\texttt l_0)
\sup_{y\in I_{\varepsilon}}\mathbb E(|\widehat b(y) - b(y)|)\leqslant
3(\texttt r_0 -\texttt l_0)\mathfrak R(\widehat b).
\end{displaymath}
{\bf Bound on $\mathbb P(\Omega_{\widehat b}^{c})$.} By Assumption \ref{assumption_b}, and by the mean value theorem,
\begin{displaymath}
b(\texttt r_{\varepsilon}) - b(\texttt l_{\varepsilon})
\leqslant -\mathfrak m_b(\texttt r_{\varepsilon} -\texttt l_{\varepsilon}),
\end{displaymath}
leading to
\begin{eqnarray*}
 \mathbb P(\Omega_{\widehat b}^{c}) & = &
 \mathbb P\left(\widehat b(\texttt r_{\varepsilon}) -\widehat b(\texttt l_{\varepsilon}) >
 -\frac{\mathfrak m_b}{2}(\texttt r_{\varepsilon} -\texttt l_{\varepsilon})\right)\\
 & \leqslant &
 \mathbb P\left(
 \widehat b(\texttt r_{\varepsilon}) -\widehat b(\texttt l_{\varepsilon}) -
 (b(\texttt r_{\varepsilon}) - b(\texttt l_{\varepsilon})) >
 \left(\mathfrak m_b -\frac{\mathfrak m_b}{2}\right)
 (\texttt r_{\varepsilon} -\texttt l_{\varepsilon})\right)\\
 & \leqslant &
 \mathbb P\left(
 |\widehat b(\texttt r_{\varepsilon}) - b(\texttt r_{\varepsilon})| +
 |\widehat b(\texttt l_{\varepsilon}) - b(\texttt l_{\varepsilon})| >
 \frac{\mathfrak m_b}{2}(\texttt r_{\varepsilon} -\texttt l_{\varepsilon})\right).
\end{eqnarray*}
Then, by Markov's inequality,
\begin{displaymath}
\mathbb P(\Omega_{\widehat b}^{c})\leqslant
\frac{2}{\mathfrak m_b(\texttt r_{\varepsilon} -\texttt l_{\varepsilon})}
\sup_{y\in I_{\varepsilon}}\mathbb E(|\widehat b(y) - b(y)|)\leqslant
\frac{2\mathfrak R(\widehat b)}{
\mathfrak m_b(\texttt r_{\varepsilon} -\texttt l_{\varepsilon})}.
\end{displaymath}
In conclusion, by Inequality (\ref{risk_bound_practical_estimator_b_1}) together with Theorem \ref{risk_bound_monotone_estimator_b},
\begin{displaymath}
\int_{I_0}\mathbb E(|\widetilde b_{\ell,h}(x) - b(x)|)dx
\lesssim
\ell + h +\mathfrak R(\widehat b) +\mathcal R(\widehat b)
\lesssim
\ell + h +\mathfrak R(\widehat b).
\end{displaymath}
%


%
\subsection{Proof of Theorem \ref{risk_bound_monotone_estimator_a}}\label{section_proof_risk_bound_monotone_estimator_a}
Since, by Assumption \ref{assumption_drift_volatility}, $|\sigma(\cdot)|\geqslant\mathfrak m_{\sigma} > 0$ on $\mathbb R$, Menozzi et al. \cite{MPZ21}, Theorem 1.2 leads to:
\begin{itemize}
 \item For every $t\in (0,T]$, the probability distribution of $X_t$ has a continuously differentiable density $f_t$ with respect to Lebesgue's measure on $\mathbb R$.
 \item $\displaystyle{\overline{\mathfrak m}_f :=\sup_{s\in [t_0,T]}\|f_s\|_{\infty} <\infty}$ and $\displaystyle{\overline{\mathfrak m}_{f'} :=\sup_{s\in [t_0,T]}\|f_s'\|_{\infty} <\infty}$.
 \item For every compact interval $I\subset\mathbb R$,
 \begin{displaymath}
 \underline{\mathfrak m}_f(I) :=
 \inf\{f_s(x)\textrm{ $;$ }(s,x)\in [t_0,T]\times I\} > 0.
 \end{displaymath}
\end{itemize}
In particular, $t\mapsto f_t(x)$ ($x\in\mathbb R$) is integrable on $[t_0,T]$, which legitimates to consider the density function $f$ defined by
\begin{displaymath}
f(x) =\frac{1}{T_0}\int_{t_0}^{T}f_s(x)ds,
\quad\forall x\in\mathbb R.
\end{displaymath}
Moreover,
\begin{itemize}
 \item For every $x\in\mathbb R$,
 \begin{equation}\label{risk_bound_monotone_estimator_a_1}
 0 < f(x)\leqslant\overline{\mathfrak m}_f.
 \end{equation}
 \item For every $x,y\in\mathbb R$,
 \begin{equation}\label{risk_bound_monotone_estimator_a_2}
 |f(y) - f(x)|\leqslant
 \frac{1}{T_0}\int_{t_0}^{T}|f_s(y) - f_s(x)|ds
 \leqslant
 \overline{\mathfrak m}_{f'}|y - x|.
 \end{equation}
 \item For every $x\in I_{2\varepsilon}$,
 \begin{equation}\label{risk_bound_monotone_estimator_a_3}
 f(x)\geqslant
 \underline{\mathfrak m}_f(I_{2\varepsilon}).
 \end{equation}
 For this reason, in the sequel, $\mathfrak m\in (0,\min\{1,\underline{\mathfrak m}_f(I_{2\varepsilon})\})$ in the definition of $\widehat a_{\eta}$.
\end{itemize}
The proof of Theorem \ref{risk_bound_monotone_estimator_a} relies on the following lemma.
%


%
\begin{lemma}\label{risk_bound_NW_estimator_a}
Under Assumption \ref{assumption_drift_volatility},
\begin{displaymath}
\sup_{x\in I_{2\varepsilon}}\mathbb E(|\widehat f_{\eta}(x) - f(x)|^2)
\lesssim
\eta^2 +\frac{1}{N\eta}
\quad\textrm{and}\quad
\sup_{x\in I_{2\varepsilon}}\mathbb E(|\widehat{af}_{\eta}(x) - (af)(x)|^2)
\lesssim
\eta^2 +\frac{1}{N\eta}.
\end{displaymath}
Then,
\begin{displaymath}
\sup_{x\in I_{2\varepsilon}}\mathbb E(|\widehat a_{\eta}(x) - a(x)|^2)
\lesssim
\eta^2 +\frac{1}{N\eta}.
\end{displaymath}
\end{lemma}
\noindent
The proof of Lemma \ref{risk_bound_NW_estimator_a} is postponed to Section \ref{section_proof_risk_bound_NW_estimator_a}.
\\
\\
By Theorem \ref{risk_bound_monotone_estimator_b} with $b = a$ and $\widehat b =\widehat a_{\eta}$,
\begin{displaymath}
\int_{I_0}\mathbb E(|\widehat a_{\ell,h,\eta}(x) - a(x)|)dx
\lesssim
\ell + h +\mathcal R(\widehat a_{\eta}),
\end{displaymath}
and by Lemma \ref{risk_bound_NW_estimator_a},
\begin{displaymath}
\mathcal R(\widehat a_{\eta})
\leqslant
(\texttt r_{2\varepsilon} -\texttt l_{2\varepsilon})
\sup_{x\in I_{2\varepsilon}}
\mathbb E(|\widehat a_{\eta}(x) - a(x)|^2)^{\frac{1}{2}}\lesssim
\left(\eta^2 +\frac{1}{N\eta}\right)^{\frac{1}{2}}.
\end{displaymath}
Thus,
\begin{displaymath}
\int_{I_0}\mathbb E(|\widehat a_{\ell,h,\eta}(x) - a(x)|)dx
\lesssim
\ell + h +\left(\eta^2 +\frac{1}{N\eta}\right)^{\frac{1}{2}}.
\end{displaymath}
%


%
\subsubsection{Proof of Lemma \ref{risk_bound_NW_estimator_a}}\label{section_proof_risk_bound_NW_estimator_a}
Consider $x\in I_{2\varepsilon}$.
\\
\\
{\bf Risk bound on $\widehat f_{\eta}(x)$.} First, by Inequality (\ref{risk_bound_monotone_estimator_a_2}),
\begin{eqnarray*}
 |\mathbb E(\widehat f_{\eta}(x)) - f(x)|^2 & = &
 \left|\frac{1}{T_0}\int_{t_0}^{T}\mathbb E(
 K_{\eta}(X_s - x))ds - f(x)\int_{-1}^{1}K(y)dy\right|^2\\
 & = &
 \left|\int_{-1}^{1}K(y)(f(\eta y + x) - f(x))dy\right|^2\leqslant
 \mathfrak c_1\eta^2
 \quad {\rm with}\quad
 \mathfrak c_1 =
 \overline{\mathfrak m}_{f'}^{2}\int_{-1}^{1}y^2K(y)dy.
\end{eqnarray*}
Now, by Inequality (\ref{risk_bound_monotone_estimator_a_1}),
\begin{eqnarray*}
 {\rm var}(\widehat f_{\eta}(x))
 & \leqslant &
 \frac{1}{NT_0}\int_{t_0}^{T}\mathbb E(K_{\eta}(X_s - x)^2)ds\\
 & = &
 \frac{1}{N}\int_{-\infty}^{\infty}K_{\eta}(y - x)^2f(y)dy\\
 & = &
 \frac{1}{N\eta}\int_{-1}^{1}K(y)^2f(\eta y + x)dy
 \leqslant
 \frac{\mathfrak c_2}{N\eta}
 \quad {\rm with}\quad
 \mathfrak c_2 =
 \overline{\mathfrak m}_f\|K\|^2.
\end{eqnarray*}
Thus,
\begin{displaymath}
\mathbb E(|\widehat f_{\eta}(x) - f(x)|^2) =
|\mathbb E(\widehat f_{\eta}(x)) - f(x)|^2 + {\rm var}(\widehat f_{\eta}(x))\leqslant
\mathfrak c_1\eta^2 +\frac{\mathfrak c_2}{N\eta}.
\end{displaymath}
Finally, since the constants $\mathfrak c_1$ and $\mathfrak c_2$ don't depend on $x$,
\begin{equation}\label{risk_bound_monotone_estimator_a_4}
\sup_{y\in I_{2\varepsilon}}\mathbb E(|\widehat f_{\eta}(y) - f(y)|^2)\leqslant
\mathfrak c_1\eta^2 +\frac{\mathfrak c_2}{N\eta}.
\end{equation}
{\bf Risk bound on $\widehat{af}_{\eta}(x)$.} First, since $a$ and $a'$ are bounded on $I_{2\varepsilon}$, and by Inequalities (\ref{risk_bound_monotone_estimator_a_1}) and (\ref{risk_bound_monotone_estimator_a_2}), $af$ is Lipschitz continuous from $I_{2\varepsilon}$ into $\mathbb R$. Then, by the martingale property of It\^o's integral with respect to the Brownian motion,
\begin{eqnarray*}
 |\mathbb E(\widehat{af}_{\eta}(x)) - (af)(x)|^2 & = &
 \left|\frac{1}{T_0}\int_{t_0}^{T}\mathbb E(
 K_{\eta}(X_s - x)a(X_s))ds - (af)(x)\int_{-1}^{1}K(y)dy\right|^2\\
 & = &
 \left|\int_{-1}^{1}K(y)((af)(\eta y + x) - (af)(x))dy\right|^2
 \leqslant\mathfrak c_3(\varepsilon)\eta^2
\end{eqnarray*}
with
\begin{displaymath}
\mathfrak c_3(\varepsilon) =
\|af' + fa'\|_{\infty,I_{2\varepsilon}}^{2}\int_{-1}^{1}y^2K(y)dy.
\end{displaymath}
Now, by the isometry property of It\^o's integral with respect to the Brownian motion, and by Inequality (\ref{risk_bound_monotone_estimator_a_1}),
\begin{eqnarray*}
 {\rm var}(\widehat{af}_{\eta}(x))
 & \leqslant &
 \frac{1}{NT_{0}^{2}}\int_{t_0}^{T}\mathbb E(K_{\eta}(X_s - x)^2\sigma(X_s)^2)ds\\
 & = &
 \frac{1}{NT_0}\int_{-\infty}^{\infty}K_{\eta}(y - x)^2\sigma(y)^2f(y)dy\\
 & = &
 \frac{1}{N\eta T_0}\int_{-1}^{1}K(y)^2(\sigma^2f)(\eta y + x)dy
 \leqslant
 \frac{\mathfrak c_4}{N\eta}
 \quad {\rm with}\quad
 \mathfrak c_4 =\frac{\|\sigma^2f\|_{\infty}\|K\|^2}{T_0}.
\end{eqnarray*}
Thus,
\begin{displaymath}
\mathbb E(|\widehat{af}_{\eta}(x) - (af)(x)|^2) =
|\mathbb E(\widehat{af}_{\eta}(x)) - (af)(x)|^2 + {\rm var}(\widehat{af}_{\eta}(x))\leqslant
\mathfrak c_3(\varepsilon)\eta^2 +\frac{\mathfrak c_4}{N\eta}.
\end{displaymath}
Finally, since the constants $\mathfrak c_3(\varepsilon)$ and $\mathfrak c_4$ don't depend on $x$,
\begin{equation}\label{risk_bound_monotone_estimator_a_5}
\sup_{y\in I_{2\varepsilon}}\mathbb E(|\widehat{af}_{\eta}(y) - (af)(y)|^2)\leqslant
\mathfrak c_3(\varepsilon)\eta^2 +\frac{\mathfrak c_4}{N\eta}.
\end{equation}
{\bf Risk bound on $\widehat a_{\eta}(x)$.} First of all,
\begin{displaymath}
\widehat a_{\eta}(x) - a(x) =
\left(\frac{\widehat{af}_{\eta}(x) - (af)(x)}{\widehat f_{\eta}(x)} +
\frac{f(x) -\widehat f_{\eta}(x)}{(\widehat f_{\eta}f)(x)}(af)(x)\right)
\mathbf 1_{\widehat f_{\eta}(x) >\frac{\mathfrak m}{2}} -
a(x)\mathbf 1_{\widehat f_{\eta}(x)\leqslant\frac{\mathfrak m}{2}}.
\end{displaymath}
Moreover, by Inequality (\ref{risk_bound_monotone_estimator_a_3}), for every $\omega\in\{\widehat f_{\eta}(x)\leqslant\frac{\mathfrak m}{2}\}$,
\begin{displaymath}
|f(x) -\widehat f_{\eta}(x,\omega)|\geqslant f(x) -\widehat f_{\eta}(x,\omega)\geqslant
\mathfrak m -\frac{\mathfrak m}{2} >\frac{\mathfrak m}{3}.
\end{displaymath}
So,
\begin{eqnarray*}
 |\widehat a_{\eta}(x) - a(x)|^2
 & \leqslant &
 \left|\frac{\widehat{af}_{\eta}(x) - (af)(x)}{\widehat f_{\eta}(x)} +
 \frac{f(x) -\widehat f_{\eta}(x)}{(\widehat f_{\eta}f)(x)}(af)(x)\right|^2
 \mathbf 1_{\widehat f_{\eta}(x) >\frac{\mathfrak m}{2}} +
 a(x)^2\mathbf 1_{|f(x) -\widehat f_{\eta}(x)| >\frac{\mathfrak m}{3}}\\
 & \leqslant &
 \frac{\mathfrak c_5(\varepsilon)}{\mathfrak m^2}(
 |\widehat{af}_{\eta}(x) - (af)(x)|^2 + |\widehat f_{\eta}(x) - f(x)|^2) +
 \|a\|_{\infty,I_{2\varepsilon}}^{2}
 \mathbf 1_{|f(x) -\widehat f_{\eta}(x)| >\frac{\mathfrak m}{3}}
\end{eqnarray*}
with
\begin{displaymath}
\mathfrak c_5(\varepsilon) =
8\max\{1,\|af\|_{\infty,I_{2\varepsilon}}^{2}\}.
\end{displaymath}
Therefore, by Markov's inequality, and by Inequalities (\ref{risk_bound_monotone_estimator_a_4}) and (\ref{risk_bound_monotone_estimator_a_5}),
\begin{eqnarray*}
 \mathbb E(|\widehat a_{\eta}(x) - a(x)|^2)
 & \leqslant &
 \frac{1}{\mathfrak m^2}(
 \mathfrak c_5(\varepsilon)\mathbb E(|\widehat{af}_{\eta}(x) - (af)(x)|^2)\\
 & &
 \hspace{1.5cm} +
 (\mathfrak c_5(\varepsilon) + 9\|a\|_{\infty,I_{2\varepsilon}}^{2})
 \mathbb E(|\widehat f_{\eta}(x) - f(x)|^2))
 \lesssim
 \eta^2 +\frac{1}{N\eta}.
\end{eqnarray*}
%


%
\subsection{Proof of Proposition \ref{risk_bound_practical_estimator_a}}\label{section_proof_risk_bound_practical_estimator_a}
By Proposition \ref{risk_bound_practical_estimator_b} with $b = a$ and $\widehat b =\widehat a_{\eta}$,
\begin{displaymath}
\int_{I_0}\mathbb E(|\widehat a_{\ell,h,\eta}(x) - a(x)|)dx
\lesssim
\ell + h +
\mathfrak R(\widehat a_{\eta}).
\end{displaymath}
Moreover, by Lemma \ref{risk_bound_NW_estimator_a},
\begin{displaymath}
\mathfrak R(\widehat a_{\eta})\lesssim
\left(\eta^2 +\frac{1}{N\eta}\right)^{\frac{1}{2}}.
\end{displaymath}
Thus,
\begin{displaymath}
\int_{I_0}\mathbb E(|\widehat a_{\ell,h,\eta}(x) - a(x)|)dx
\lesssim
\ell + h +\left(\eta^2 +\frac{1}{N\eta}\right)^{\frac{1}{2}}.
\end{displaymath}
%


%

%
\end{document}